\input amstex
\input amsppt.sty
\magnification=\magstep1
\hsize=30truecc
\vsize=22.2truecm
\baselineskip=16truept
\NoBlackBoxes
\TagsOnRight \pageno=1 \nologo
\def\Z{\Bbb Z}
\def\N{\Bbb N}

\def\l{\left}
\def\r{\right}
\def\bg{\bigg}
\def\({\bg(}
\def\[{\bg\lfloor}
\def\){\bg)}
\def\]{\bg\rfloor}
\def\t{\text}
\def\f{\frac}

\def\se {\subseteq}

\def\sm{\setminus}

\def\bi{\binom}
\def\eq{\equiv}

\def\ls{\leqslant}
\def\gs{\geqslant}
\def\mo{\roman{mod}}

\def\al{\alpha}
\def\da{\delta}

\def\Proof{\noindent{\it Proof}}

\def\Ack{\medskip\noindent {\bf Acknowledgment}}
\hbox {Acta Arith. 148(2011), 63-76.}
\bigskip
\topmatter
\title $p$-adic valuations of some sums of multinomial coefficients\endtitle
\author Zhi-Wei Sun\endauthor
\leftheadtext{Zhi-Wei Sun} \rightheadtext{$p$-adic valuations of
sums of multinomial coefficients}
\affil Department of Mathematics, Nanjing University\\
 Nanjing 210093, People's Republic of China
  \\  zwsun\@nju.edu.cn
  \\ {\tt http://math.nju.edu.cn/$\sim$zwsun}
\endaffil
\abstract Let $m$ and $n>0$ be integers. Suppose that $p$
is an odd prime dividing $m-4$. We show that
$$\nu_p\(\sum_{k=0}^{n-1}\f{\bi{2k}k}{m^k}\)\gs\nu_p(n)\ \ \t{and}
\ \ \nu_p\(\sum_{k=0}^{n-1}\bi{n-1}k(-1)^k\f{\bi{2k}k}{m^k}\)\gs\nu_p(n),$$
where $\nu_p(x)$ denotes the $p$-adic valuation of $x$.
Furthermore, if $p>3$ then
$$\f1n\sum_{k=0}^{n-1}\f{\bi{2k}k}{m^k}\eq\f{\bi{2n-1}{n-1}}{4^{n-1}}\ (\mo\ p^{\nu_p(m-4)})$$
and
$$\f1n\sum_{k=0}^{n-1}\bi{n-1}k(-1)^k\f{\bi{2k}k}{m^k}\eq\f{C_{n-1}}{4^{n-1}}\ (\mo\ p^{\nu_p(m-4)}),$$
where $C_k$ denotes the Catalan number $\f1{k+1}\bi{2k}k$.
This implies several conjectures of Guo and Zeng [GZ].
We also raise two conjectures, and prove that $n>1$ is a prime if and only if
$$\sum_{k=0}^{n-1}\bi{(n-1)k}{k,\ldots,k}\eq0\ (\mo\ n),$$
where $\bi{k_1+\cdots+k_{n-1}}{k_1,\ldots,k_{n-1}}$ denotes the multinomial coefficient
$\f{(k_1+\cdots+k_{n-1})!}{k_1!\cdots k_{n-1}!}$.
\endabstract
\thanks 2010 {\it Mathematics Subject Classification}.\,Primary 11B65;
Secondary 05A10,\,11A07, 11S99.
\newline\indent Supported by the National Natural Science
Foundation (grant 10871087) and the Overseas Cooperation Fund (grant 10928101) of China.
\endthanks
\endtopmatter
\document

\heading{1. Introduction}\endheading

Let $p$ be a prime. In 2006 Pan and Sun [PS]
obtained various congruences modulo $p$ involving central binomial coefficients and Catalan numbers.
Later Sun and Tauraso [ST1, ST2] made some further refinements;
for example, they proved that for any $a\in\Z^+=\{1,2,3,\ldots\}$ we have
$$\sum_{k=0}^{p^a-1}\bi{2k}k\eq\l(\f{p^a}3\r)\ (\mo\ p^2).$$
Recently the author [S10] managed to determine
$\sum_{k=0}^{p^a-1}\bi{2k}{k}/m^k$ mod $p^2$ for any integer $m$
not divisible by $p$.

 Motivated by the above work, Guo and Zeng [GZ] obtained some congruences involving central $q$-binomial coefficients
 and raised several conjectures on $p$-adic valuations of some sums of binomial coefficients.

Throughout the paper, for a prime $p$,  the $p$-adic valuation (or $p$-adic order)
of an integer $m$ is given by
$$\nu_p(m)=\sup\{a\in\Z:\ p^a\mid m\},$$
and we define
$\nu_p(m/n)=\nu_p(m)-\nu_p(n)$ for any $m\in\Z$ and $n\in\Z^+$.
For example,
$$\nu_2\l(\f23\r)=\nu_2(2)-\nu_2(3)=1
\ \ \t{and}\ \ \nu_3\l(\f49\r)=\nu_3(4)-\nu_3(9)=-2.$$
For an assertion $A$ we adopt the Iverson notation:
$$[A]=\cases 1&\t{if}\ A\ \t{holds},\\0&\t{otherwise}.\endcases$$
Thus $[m=n]$ coincides with the Kronecker symbol $\da_{m,n}$.

Our following result implies several conjectures of Guo and Zeng [GZ, Section 5].

\proclaim{Theorem 1.1} Let $m\in\Z$ and $n\in\Z^+$. Suppose that $p$
is an odd prime dividing $m-4$. Then
$$\nu_p\(\sum_{k=0}^{n-1}\f{\bi{2k}k}{m^k}\)\gs \nu_p(n)
\ \t{and}\ \nu_p\(\sum_{k=0}^{n-1}\bi{n-1}k(-1)^k\f{\bi{2k}k}{m^k}\)\gs \nu_p(n).\tag1.1$$
Furthermore,
$$\f1n\sum_{k=0}^{n-1}\f{\bi{2k}k}{m^k}\eq\f{\bi{2n-1}{n-1}}{4^{n-1}}
+\da_{p,3}[3\mid n]\f{m-4}3\bi{2n/3^{\nu_3(n)}-1}{n/3^{\nu_3(n)}-1}\ (\mo\ p^{\nu_p(m-4)})\tag1.2$$
and also
$$\f1n\sum_{k=0}^{n-1}\bi{n-1}k(-1)^k\f{\bi{2k}k}{m^k}\eq\f{C_{n-1}}{4^{n-1}}\ (\mo\ p^{\nu_p(m-4)-\da_{p,3}}),\tag1.3$$
where $C_k$ denotes the Catalan number $\f1{k+1}\bi{2k}k=\bi{2k}k-\bi{2k}{k+1}$.
Thus, for $a\in\Z^+$ we have
$$\f1{p^a}\sum_{k=0}^{p^a-1}\f{\bi{2k}k}{m^k}\eq 1+\da_{p,3}\f{m-4}3\eq\f{m-1}3\ (\mo\ p)\tag1.4$$
and also
$$\f1{p^a}\sum_{k=0}^{p^a-1}\bi{p^a-1}k(-1)^k\f{\bi{2k}k}{m^k}\eq-1\ (\mo\ p)\quad \t{provided}\ p\not=3.\tag1.5$$
\endproclaim

Now we give various consequences of Theorem 1.1.

\proclaim{Corollary 1.1 {\rm ([GZ, Conjecture 5.1])}} Let $p$ be a
prime divisor of $4m-1$ with $m\in\Z$. Then
$$\nu_p\(\sum_{k=0}^{n-1}\bi{2k}k m^k\)\gs\nu_p(n)\tag1.6$$
for all $n\in\Z^+$.
\endproclaim
\Proof. As $p\nmid m$, there exists an integer $m_*$ such that $m_*m\eq1\ (\mo\ p^{\nu_p(n)})$
and hence $m_*\eq4\ (\mo\ p)$. By Theorem 1.1, for any $n\in\Z^+$ we have
$$\sum_{k=0}^{n-1}\bi{2k}km^k\eq\sum_{k=0}^{n-1}\f{\bi{2k}k}{m_*^k}\eq0\ (\mo\ p^{\nu_p(n)}).$$
This concludes the proof. \qed

\proclaim{Corollary 1.2 {\rm ([GZ, Conjecture 5.2])}} Let $n=|4m-1|$
with $m\in\Z$. Then
$$\sum_{k=0}^{n-1}\bi{2k}k m^k\eq0\ (\mo\ n).\tag1.7$$
\endproclaim
\Proof. By Corollary 1.1, (1.6) holds for any prime $p$ dividing $n$. So (1.7) is valid. \qed

\proclaim{Corollary 1.3 {\rm ([GZ, Conjecture 5.4])}} Let $p>3$ be a prime and $a\in\Z^+$. Then
$$\sum_{k=0}^{p^a-1}\bi{2k}k\l(\f{1-(-1)^{(p-1)/2}p}4\r)^k\eq p^a\ (\mo\ p^{a+1}).\tag1.8$$
\endproclaim
\Proof. Let $m=(1-(-1)^{(p-1)/2}p)/4$. Then $m\in\Z$ and $p\mid 4m-1$.
Choose an integer $m_*$ such that $mm_*\eq1\ (\mo\ p^{a+1})$. Clearly $m^*\eq 4\ (\mo\ p)$.
Applying Theorem 1.1 we get
$$\f1{p^a}\sum_{k=0}^{p^a-1}\bi{2k}k m^k\eq\f1{p^a}\sum_{k=0}^{p^a-1}\f{\bi{2k}k}{m_*^k}\eq1
\ (\mo\ p).$$
So (1.8) holds. \qed

Note that (1.8) in the case $p=5$ yields
$$\sum_{k=0}^{5^a-1}(-1)^k\bi{2k}k\eq 5^a\ (\mo\ 5^{a+1}),\tag1.9$$
which is the second congruence in [GZ, Conjecture 3.5].
\medskip

\proclaim{Corollary 1.4 {\rm ([GZ, Conjecture 5.3])}} For $a\in\Z^+$ we have
$$\align \sum_{k=0}^{3^a-1}(-2)^k\bi{2k}k\eq& 3^a\ (\mo\ 3^{a+1}),\tag1.10
\\\sum_{k=0}^{3^a-1}(-5)^k\bi{2k}k\eq&-3^a\ (\mo\ 3^{a+1}),\tag1.11
\\\sum_{k=0}^{7^a-1}(-5)^k\bi{2k}k\eq& 7^a\ (\mo\ 7^{a+1}).\tag1.12
\endalign$$
\endproclaim
\Proof. Choose integers $m_1,m_2,m_3$ such that
$$m_1\eq-\f12\ (\mo\ 3^{a+1}),\ m_2\eq-\f15\ (\mo\ 3^{a+1}),\ m_3\eq -\f15\ (\mo\ 7^{a+1}).$$
Then
$$m_1\eq 4\ (\mo\ 3^2),\ \ m_2\eq 4\ (\mo\ 3)\ \ \t{and}\ \ m_3\eq 4\ (\mo\ 7).$$
So it suffices to apply (1.4). \qed

 \medskip

 (1.4) in the case $p=3$, together with our computation via {\tt Mathematica}, leads us to raise the following conjecture.

 \proclaim{Conjecture 1.1} Let $m\in\Z$ with $m\eq1\ (\mo\ 3)$. Then
 $$\nu_3\(\f1n\sum_{k=0}^{n-1}\f{\bi{2k}k}{m^k}\)\gs \min\{\nu_3(n),\nu_3(m-1)-1\}\tag1.13$$
 and
$$\nu_3\(\f1n\sum_{k=0}^{n-1}\bi{n-1}k(-1)^k\f{\bi{2k}k}{m^k}\)\gs \min\{\nu_3(n),\nu_3(m-1)\}-1\tag1.14$$
for every $n\in\Z^+$.  Furthermore,
 $$\f1{3^a}\sum_{k=0}^{3^a-1}\f{\bi{2k}k}{m^k}\eq\f{m-1}3\ (\mo\ 3^{\nu_3(m-1)})$$
 for any integer $a\gs\nu_3(m-1)$, and
 $$\f1{3^a}\sum_{k=0}^{3^a-1}\bi{3^a-1}k(-1)^k\f{\bi{2k}k}{m^k}\eq-\f{m-1}3\ (\mo\ 3^{\nu_3(m-1)})$$
for each integer $a>\nu_3(m-1)$. Also,
$$\sum_{k=0}^{3^a-1}\bi{3^a-1}k(-1)^k\bi{2k}k\eq -3^{2a-1}\ (\mo\ 3^{2a})\ \ \t{for every}\ a=2,3,\ldots.$$
 \endproclaim

 We remark that Strauss, Shallit and Zagier [SSZ] used a special technique to
 show that for any $n\in\Z^+$ we have
$$\nu_3\(\sum_{k=0}^{n-1}\bi{2k}k\)=2\nu_3(n)+\nu_3\l(\bi{2n}n\r).$$

For any $k\in\N=\{0,1,2,\ldots\}$, the central binomial coefficient $\bi{2k}k$ coincides with the multinomial coefficient $\bi{2k}{k,k}$.
In general, the multinomial coefficient
$$\bi{k_1+\cdots+k_n}{k_1,\ldots,k_n}=\f{(k_1+\cdots+k_{n})!}{k_1!\cdots k_{n}!}$$
in the case $k_1,\ldots,k_n=k\in\N$ gives
$$\bi{nk}{k,\ldots,k}=\f{(nk)!}{(k!)^n}.$$

 Now we pose one more conjecture which involves multinomial coefficients.

 \proclaim{Conjecture 1.2} For any prime $p$ and positive integer $n$ we have
 $$\nu_p\(\sum_{k=0}^{n-1}\bi{(p-1)k}{k,\ldots,k}\)\gs\nu_p(n)\tag1.15$$
 and
 $$\nu_p\(\sum_{k=0}^{n-1}\bi{n-1}k(-1)^k\bi{(p-1)k}{k,\ldots,k}\)\gs\nu_p(n).\tag1.16$$
 Furthermore, $\nu_p(n)$ in $(1.15)$ can be replaced by
 $\nu_p(n\bi{2n}n)$ if $p>2$.
 \endproclaim

Observe that $$\f{(4k)!}{(k!)^4}=\bi{4k}{2k}\bi{2k}k^2$$ and hence
(1.15) in the case $p=5$ yields the first congruence in [GZ,
Conjecture 5.6].

Concerning Conjecture 1.2 we can prove the following result.

\proclaim{Theorem 1.2}  Let $p$ be a prime.

{\rm (i)} We have
$$\sum_{k=0}^{p-1}\bi{(p-1)k}{k,\ldots,k}\eq pB_{p-1}+(-1)^{p-1}-2p\ (\mo\ p^2),\tag1.17$$
where $B_n$ denotes the $n$th Bernoulli number. Also, an integer $m>1$ is a prime if and only if
$$\sum_{k=0}^{m-1}\bi{(m-1)k}{k,\ldots,k}\eq0\ (\mo\ m).\tag1.18$$

{\rm (ii)}  Let $n\in\Z^+$.
If $n\not\eq1\ \mo\ p)$ or there is a digit greater than $1$ in the representation of $n$
in base $p$, then
$$\sum_{k=0}^{n-1}\bi{(p-1)k}{k,\ldots,k}\eq0\ (\mo\ p),\tag1.19$$
otherwise we have
$$\sum_{k=0}^{n-1}\bi{(p-1)k}{k,\ldots,k}\eq(-1)^{\psi_p(n)-1}\ (\mo\ p),\tag1.20$$
where $\psi_p(n)$ denotes the sum of all the digits in the representation of $n$ in base $p$.

{\rm (iii)}  $(1.15)$ holds for all $n\in\Z^+$ if and only if so does $(1.16)$.
\endproclaim

A basic problem in number theory is to characterize primes.
However, besides the well-known Wilson theorem, no other simple
congruence characterization of primes has been proved before. Thus
our characterization of primes via (1.18) is particularly interesting.

It is curious to know what odd primes $p$ satisfy the congruence
$$\sum_{k=0}^{p-1}\bi{(p-1)k}{k,\ldots,k}\eq0\ (\mo\ p^2)\ \ (\t{i.e.},\ pB_{p-1}\eq2p-1\ (\mo\ p^2)).$$
Using {\tt Mathematica} we only find four such primes (they are 3, 11, 107, 4931) among the first 15,000 primes.
It seems that all such primes are congruent to 3 modulo 8. From the proof of (1.17) we see that such odd primes
are exactly those odd primes $p$ satisfying $(p-2)!\eq1\ (\mo\ p^2)$, which were investigated by P. Saridis [S]
who also found the above four primes. (The author thanks Prof. N.J.A. Sloane for informing him about the reference [S].)

In the next section we are going to provide some lemmas.
Theorems 1.1 and 1.2 will be proved in Sections 3 and 4 respectively.

\heading{2. Some Lemmas}\endheading

\proclaim{Lemma 2.1 {\rm ([ST1, Theorem 2.1])}} For any $n\in\Z^+$
and $d\in\Z$, we have
$$\aligned&\sum_{0\ls k<n}\bi{2k}{k+d}x^{n-1-k}+[d>0]x^nu_d(x-2)
\\&\ \ =\sum_{0\ls k<n+d}\bi{2n}ku_{n+d-k}(x-2),
\endaligned\tag2.1$$
where the polynomial sequence $\{u_k(x)\}_{k\gs0}$ is defined as
follows:
$$u_0(x)=0,\ u_1(x)=1,\ \t{and}\ u_{k+1}(x)=xu_k(x)-u_{k-1}(x)\
(k=1,2,3,\ldots).$$
\endproclaim

Let $A,B\in\Z$. The Lucas sequence $u_n=u_n(A,B)\ (n\in\N)$
is defined by
$$u_0=0,\ u_1=1,\ \t{and}\ u_{n+1}=Au_n-Bu_{n-1}\ (n=1,2,3,\ldots).$$
The characteristic equation $x^2-Ax+B=0$ has two roots
$$\al=\f{A+\sqrt{\Delta}}2\quad\t{and}\quad\beta=\f{A-\sqrt{\Delta}}2,$$
where $\Delta=A^2-4B$. It is well known that for any $n\in\N$ we have
$$u_n=\sum_{0\ls k<n}\al^k\beta^{n-1-k}\quad\t{and hence}\quad (\al-\beta)u_n=\al^n-\beta^n.$$
The reader may consult [S06] for connections between Lucas sequences and quadratic fields.

\proclaim{Lemma 2.2} Let $A,B\in\Z$ and let $d\in\Z^+$ be an odd divisor of $\Delta=A^2-4B$.
Then, for any $n\in\Z^+$, we have
$$\f{u_n(A,B)}n\eq\l(\f A2\r)^{n-1}+\cases(A/2)^{n-3}\Delta/3\ (\mo\ d)&\t{if}\ 3\mid d\ \t{and}\ 3\mid n,
\\0\ (\mo\ d)&\t{otherwise}.\endcases\tag2.2$$
\endproclaim
\Proof. When $\Delta=0$, by induction $u_k(A,B)=k(A/2)^{k-1}$ for
all $k\in\Z^+$, and hence the desired result follows.

Now we assume that $\Delta\not=0$. Then
$$\align u_n(A,B)=&\f1{\sqrt{\Delta}}\(\l(\f{A+\sqrt{\Delta}}2\r)^n-\l(\f{A-\sqrt{\Delta}}2\r)^n\)
\\=&\f2{2^n}\sum\Sb 0\ls k\ls n\\2\nmid k\endSb\bi nkA^{n-k}\Delta^{(k-1)/2}
\\=&\f1{2^{n-1}}\sum\Sb 1\ls k\ls n\\2\nmid k\endSb\f nk\bi {n-1}{k-1}A^{n-k}\Delta^{(k-1)/2}
\endalign$$
and hence
$$\f{u_n(A,B)}n-\l(\f A2\r)^{n-1}
=\sum\Sb 1<k\ls n\\2\nmid k\endSb\bi{n-1}{k-1}\l(\f A2\r)^{n-k}\f{\Delta^{(k-1)/2}}{k2^{k-1}}.\tag2.3$$

 For $k=5,7,9,\ldots$, clearly
$k<3^{(k-1)/2}$ and hence $\nu_p(k)\ls(k-3)/2$ for any prime divisor $p$ of $d$,
thus $\Delta\Delta^{(k-3)/2}/k\eq0\ (\mo\ d)$.
Note also that
$$\align &\bi{n-1}{3-1}\l(\f A2\r)^{n-3}\f{\Delta^{(3-1)/2}}{3\times 2^{3-1}}
\\=&\f{(n-1)(n-2)}2\l(\f A2\r)^{n-3}\f{\Delta}{3\times 4}
\\\eq&\cases(A/2)^{n-3}\Delta/3\ (\mo\ d)&\t{if}\ 3\mid d\ \t{and}\ 3\mid n,
\\0\ (\mo\ d)&\t{otherwise}.\endcases
\endalign$$
So (2.2) follows from (2.3).

The proof of Lemma 2.2 is now complete. \qed

\proclaim{Lemma 2.3} If $p$ is a prime, and
$$a=\sum_{i=0}^ka_ip^i\ \ \t{and}\ \ b=\sum_{i=0}^kb_ip^i\ \ (a_i,b_i\in\{0,\ldots,p-1\}),$$
then we have the Lucas congruence
$$\bi ab\eq\prod_{i=0}^k\bi{a_i}{b_i}\ (\mo\ p).$$
\endproclaim

This lemma is a well-known result due to Lucas, see, e.g., [St, p.\,44].

\proclaim{Lemma 2.4} Let $p$ be a prime and let $h\in\Z^+$ and $m\in\Z\sm\{0\}$. Then we have
$$\min_{1\ls k\ls n}\nu_p\(\f1k\sum_{l=0}^{k-1}\bi{k-1}l(-1)^l\f{\bi{hl}{l,\ldots,l}}{m^l}\)
=\min_{1\ls k\ls n}\nu_p\(\f1k\sum_{l=0}^{k-1}\f{\bi{hl}{l,\ldots,l}}{m^l}\)\tag2.4$$
for every $n=1,2,3,\ldots$.
\endproclaim
\Proof. By a confirmed conjecture of Dyson (cf. [D, Go, Z] or [St, p.\,44]), for any $k\in\N$ the constant term
of the Laurent polynomial $$\prod\Sb 1\ls i,j\ls h\\i\not=j\endSb\l(1-\f{x_i}{x_j}\r)^k$$
coincides with the multinomial coefficient $\bi{hk}{k,\ldots,k}$.

Let $n\in\Z^+$. Then
$$\align&\sum_{k=0}^{n-1}\f1{m^k}\prod\Sb 1\ls i,j\ls h\\i\not=j\endSb\l(1-\f{x_i}{x_j}\r)^k
\\=&\f{(m^{-1}\prod_{1\ls i,j\ls h,\, i\not=j}(1-x_i/x_j))^n-1}{m^{-1}\prod_{1\ls i,j\ls h,\, i\not=j}(1-x_i/x_j)-1}
\\=&\sum_{k=1}^{n}\bi nk\(\f1m\prod\Sb 1\ls i,j\ls h\\i\not=j\endSb\l(1-\f{x_i}{x_j}\r)-1\)^{k-1}
\\=&\sum_{k=1}^n\f nk\bi{n-1}{k-1}\sum_{l=0}^{k-1}\bi{k-1}l\f{(-1)^{k-1-l}}{m^l}\prod\Sb1\ls i,j\ls h\\i\not=j\endSb
\l(1-\f{x_i}{x_j}\r)^l.
\endalign$$
Comparing the constant terms of both sides we get
$$\f1n\sum_{k=0}^{n-1}\f{\bi{hk}{k,\ldots,k}}{m^k}
=\sum_{k=1}^n\bi{n-1}{k-1}\f {(-1)^{k-1}}k\sum_{l=0}^{k-1}\bi{k-1}l(-1)^{l}\f{\bi{hl}{l,\ldots,l}}{m^l}.\tag2.5$$

Recall that for any sequences $\{a_n\}_{n\gs0}$ and $\{b_n\}_{n\gs0}$ of complex numbers we have
$$\align &a_n=\sum_{k=0}^n\bi nk(-1)^kb_k\quad\t{for all}\ n=0,1,2,\ldots
\\\iff&b_n=\sum_{k=0}^n\bi nk(-1)^ka_k\quad\t{for all}\ n=0,1,2,\ldots.
\endalign$$
(See, e.g., [R, p.\,43].)
So (2.5) holds for all $n\in\Z^+$ if and only if for each $n\in\Z^+$ we have
$$\sum_{k=1}^n\bi{n-1}{k-1}\f{(-1)^{k-1}}k\sum_{l=0}^{k-1}\f{\bi{hl}{l,\ldots,l}}{m^l}
=\f1n\sum_{l=0}^{n-1}\bi{n-1}l(-1)^l\f{\bi{hl}{l,\ldots,l}}{m^l}.
\tag2.6$$

Since both (2.5) and (2.6) are valid for all $n\in\Z^+$, (2.4) holds for any $n\in\Z^+$.
This concludes the proof. \qed

\heading{3. Proof of Theorem 1.1}\endheading

 Observe that $p\nmid m$ since $p\mid m-4$ and $p\not=2$.
Applying Lemma 2.1 with $x=m$ and $d=0$,
we get
$$\align \f{m^{n-1}}n\sum_{k=0}^{n-1}\f{\bi{2k}k}{m^k}
=&\f1n\sum_{k=0}^{n-1}\bi{2n}ku_{n-k}(m-2,1)
\\=&\sum_{k=0}^{n-1}\(2\bi{2n-1}k-\bi{2n}k\)\f{u_{n-k}(m-2,1)}{n-k}.
\endalign$$
Since $m-2\eq 2\ (\mo\ p^{\nu_p(m-4)})$, we have
$$\sum_{k=0}^{n-1}\(2\bi{2n-1}k-\bi{2n}k\)\l(\f{m-2}2\r)^{n-k-1}\eq\Sigma\ (\mo\ p^{\nu_p(m-4)})$$
where
$$\Sigma:=\sum_{k=0}^{n-1}\(2\bi{2n-1}k-\bi{2n}k\)=\bi{2n-1}{n-1}.$$
Thus, by Lemma 2.2 and the above,
$$\align&\f{m^{n-1}}n\sum_{k=0}^{n-1}\f{\bi{2k}k}{m^k}-\bi{2n-1}{n-1}
\\\eq&\da_{p,3}\sum\Sb k=0\\3\mid n-k\endSb^{n-1}
\(2\bi{2n-1}k-\bi{2n}k\)\l(\f{m-2}2\r)^{(n-k)-3}\f{m(m-4)}3
\\\eq&\da_{p,3}\f{m-4}3 S_n \ (\mo\ p^{\nu_p(m-4)})\ \ (\t{since}\ m\eq 4\ (\mo\ p^{\nu_p(m-4)})),
\endalign$$
where
$$S_n=\sum\Sb k=0\\3\mid n-k\endSb^{n-1}
\(2\bi{2n-1}k-\bi{2n}k\).$$

In the case $3\nmid n$,  for any $k\in\{0,\ldots,n-1\}$ with $k\eq n\ (\mo\ 3)$ we have
$$2\bi{2n-1}k-\bi{2n}k=\f{n-k}n\bi{2n}k\eq0\ (\mo\ 3).$$
So $3\mid S_n$ if $3\nmid n$.

In the case $3\mid n$, by Lemma 2.3, for $k\in\N$ we have
$$\bi{2n}{3k}\eq\bi{2n/3}{k}\ (\mo\ 3)$$
and
$$\align\bi{2n-1}{3k}=&\f{(2n-1)(2n-2)}{(2n-3k-1)(2n-3k-2)}\bi{2n-3}{3k}
\\\eq&\bi{2n-3}{3k}\eq\bi{2n/3-1}k\ (\mo\ 3),
\endalign$$
thus
$$\align S_n=&\sum_{k=0}^{n/3-1}\(2\bi{2n-1}{3k}-\bi{2n}{3k}\)
\\\eq&-\sum_{k=0}^{n/3-1}\(\bi{2n/3-1}{k}+\bi{2n/3}{k}\) \ (\mo\ 3)
\endalign$$
and hence
$$S_n\eq -2^{2n/3-2}-2^{2n/3-1}+\f12\bi{2n/3}{n/3}\eq\f12\bi{2q}q=\bi{2q-1}{q-1}\ (\mo\ 3)$$
with $q=n/3^{\nu_3(n)}$.

Combining the above we get
$$\align\f1n\sum_{k=0}^{n-1}\f{\bi{2k}k}{m^k}
\eq&\f{\bi{2n-1}{n-1}+\da_{p,3}[3\mid n]\f{m-4}3\bi{2q-1}{q-1}}{m^{n-1}}
\\\eq&\f{\bi{2n-1}{n-1}}{4^{n-1}}+\da_{p,3}[3\mid n]\f{m-4}3\bi{2q-1}{q-1}\ (\mo\ p^{\nu_p(m-4)}).
\endalign$$
This, together with (2.6) in the case $h=2$, yields
$$\f1n\sum_{k=0}^{n-1}\bi{n-1}k(-1)^k\f{\bi{2k}k}{m^k}
=\sigma\ (\mo\ p^{\nu_p(m-4)-\da_{p,3}}),$$
where
$$\align \sigma:=&\sum_{k=1}^n\bi{n-1}{k-1}\f{(-1)^{k-1}}{4^{k-1}}\bi{2k-1}{k-1}
\\=&-2\sum_{k=0}^n\bi{n-1}{n-k}\bi{-1/2}k=-2\bi{n-3/2}n=\f{C_{n-1}}{4^{n-1}}
\endalign$$
with the help of the Chu-Vandermonde identity (see (5.22) of [GKP, p.\,169]).

Clearly, if $n=p^a$ for some $a\in\Z^+$ then
$$\f{\bi{2n-1}{n-1}}{4^{n-1}}\eq \bi{2p^a-1}{p^a-1}=\prod_{k=1}^{p^a-1}\l(1+\f{p^a}k\r)\eq1\ (\mo\ p)$$
and
$$\f{C_{n-1}}{4^{n-1}}\eq \f1{p^a}\bi{2p^a-2}{p^a-1}=\f1{2p^a-1}\bi{2p^a-1}{p^a}\eq-1\ (\mo\ p).$$
This concludes our proof of Theorem 1.1.

\heading{4. Proof of Theorem 1.2}\endheading

\medskip

\proclaim{Lemma 4.1} Let $p$ be a prime and let $n\in\Z^+$.
If all the digits in the representation of $n$ in base $p$ belong to $\{0,1\}$, then
$$\prod_{j=1}^{p-1}\bi {jn}n\eq (-1)^{\psi_p(n)}\ (\mo\ p)$$
(where $\psi_p(n)$ is defined as in Theorem 1.2),
otherwise we have
$$\prod_{j=1}^{p-1}\bi {jn}n\eq0\ (\mo\ p).$$
\endproclaim
\Proof. Suppose that $n=\sum_{i=0}^ka_ip^i$ with $a_0,\ldots,a_k\in\{0,\ldots,p-1\}$.

If $a_0,\ldots,a_k\in\{0,1\}$ then $ja_i\ls j<p$ for all $i=0,\ldots,k$ and  $j=1,\ldots,p-1$, thus
$$\align&\prod_{j=1}^{p-1}\bi{jn}n=\prod_{j=1}^{p-1}\bi{\sum_{i=0}^k(ja_i)p^i}{\sum_{i=0}^k a_ip^i}
\\\eq&\prod_{j=1}^{p-1}\prod_{i=0}^k\bi{ja_i}{a_i}=\prod_{i=0}^k\prod_{j=1}^{p-1}\bi{ja_i}{a_i}\ \ (\t{by Lemma 2.3})
\\\eq&((p-1)!)^{|\{0\ls i\ls k:\, a_i=1\}|}\eq(-1)^{\psi_p(n)}\ (\mo\ p)\ (\t{by Wilson's theorem)}.
\endalign$$

Now assume that $\{a_0,\ldots,a_k\}\not\se\{0,1\}$. We want to show that $p\mid\bi{jn}n$ for some $j\in\{1,\ldots,p-1\}$.
Set $s=\min\{0\ls i\ls k:\, a_i>1\}$.
As $1<a_s<p$, we may choose $j\in\{1,\ldots,p-1\}$ such that $ja_s\eq1\ (\mo\ p)$. Thus
$$jn=\sum_{s<i\ls k}(ja_i)p^i+(ja_s-1)p^{s}+p^s+\sum_{0\ls t<s}(ja_t)p^t.$$
Write $$\sum_{s<i\ls k}(ja_i)p^i+(ja_s-1)p^{s}=\sum_{s<i\ls k}b_ip^i+bp^{k+1}$$
with $b_i\in\{0,\ldots,p-1\}$ and $b\in\N$. Then, with the help of Lemma 2.3, we have
$$\align\bi{jn}n=&\bi{bp^{k+1}+\sum_{s<i\ls k}b_ip^i+p^s+\sum_{0\ls t<s}(ja_t)p^t}{\sum_{i=0}^ka_ip^i}
\\\eq&\prod_{s<i\ls k}\bi{b_i}{a_i}\times\bi1{a_s}\times\prod_{0\ls t<s}\bi{ja_t}{a_t}=0\ (\mo\ p).
\endalign$$

Combining the above we have proved the desired result. \qed

\medskip

\noindent{\it Proof of Theorem 1.2}. (i) If $n$ is an integer greater than $1$, then $(pn-1)!\eq0\ (\mo\ p)$ and hence
$$\align\sum_{k=0}^{pn-1}\bi{(pn-1)k}{k,\ldots,k}=&\sum_{k=0}^{pn-1}\prod_{j=1}^{pn-1}\bi{jk}k
=1+\sum_{k=1}^{pn-1}\prod_{j=1}^{pn-1}\l(j\bi{jk-1}{k-1}\r)
\\=&1+(pn-1)!\sum_{k=1}^{pn-1}\prod_{j=1}^{pn-1}\bi{jk-1}{k-1}
\eq1\ (\mo\ p).
\endalign$$
So (1.18) fails for any composite number $m>1$.

If $1<k\ls p-1$, then
$(p-1)k\gs2(p-1)\gs p$ and hence
$$\bi{(p-1)k}{k,\ldots,k}=\f{((p-1)k)!}{(k!)^{p-1}}\eq0\ (\mo\ p).$$
Thus
$$\sum_{k=0}^{p-1}\bi{(p-1)k}{k,\ldots,k}\eq\sum_{k=0}^1\bi{(p-1)k}{k,\ldots,k}=1+(p-1)!\eq0\ (\mo\ p)$$
with the help of Wilson's theorem.

 Now we determine $\sum_{k=0}^{p-1}\bi{(p-1)k}{k,\ldots,k}$ modulo $p^2$.

 In the case $p=2$, as $B_1=-1/2$ we have
$$\sum_{k=0}^{p-1}\bi{(p-1)k}{k,\ldots,k}=1+(p-1)!=2\eq 2B_{p-1}+(-1)^{p-1}-2p\ (\mo\ p^2).$$

Now let $p$ be an odd prime. If $2<k\ls p-1$, then there exist $j_1,j_2\in\{1,\ldots,p-1\}$ such that
$j_1k\eq1\ (\mo\ p)$ and $j_2k\eq 2\ (\mo\ p)$, hence
$\bi{j_1k}k\eq\bi{j_2k}k\eq0\ (\mo\ p)$
by Lemma 2.3, and thus
$$\bi{(p-1)k}{k,\ldots,k}=\prod_{j=1}^{p-1}\bi{jk}k\eq0\ (\mo\ p^2).$$
Note also that
$$\bi{(p-1)2}{2,\ldots,2}=\prod_{j=1}^{p-1}\bi{2j}2=\prod_{j=1}^{p-1}\l(j(2j-1)\r)
\eq p!(p-2)!\eq-p\ (\mo\ p^2).$$
Therefore
$$\sum_{k=0}^{p-1}\bi{(p-1)k}{k,\ldots,k}\eq\sum_{k=0}^1\bi{(p-1)k}{k,\ldots,k}-p\eq 1+(p-1)!-p\ (\mo\ p^2)$$
and hence we have (1.17) with the help of Glaisher's result $(p-1)!\eq pB_{p-1}-p\ (\mo\ p^2)$ (cf. [Gl]).

(ii) Write $n=pm+r$ with $m\in\N$ and $r\in\{0,\ldots,p-1\}$. If $m>0$ then
$$\align\sum_{k=0}^{pm-1}\bi{(p-1)k}{k,\ldots,k}=&\sum_{k=0}^{pm-1}\prod_{j=1}^{p-1}\bi{jk}k
=\sum_{k=0}^{m-1}\sum_{t=0}^{p-1}\prod_{j=1}^{p-1}\bi{pjk+jt}{pk+t}
\\\eq&\sum_{k=0}^{m-1}\sum_{t=0}^1\prod_{j=1}^{p-1}\bi{pjk+jt}{pk+t}\ (\t{by Lemma 4.1})
\\\eq&\sum_{k=0}^{m-1}\sum_{t=0}^1\prod_{j=1}^{p-1}\(\bi{jt}t\bi{jk}k\)\ (\t{by Lemma 2.3})
\\\eq&\sum_{k=0}^{m-1}(1+(p-1)!)\prod_{j=1}^{p-1}\bi{jk}k\eq0\ (\mo\ p).
\endalign$$
Similarly,
$$\sum_{pm\ls k<pm+r}\bi{(p-1)k}{k,\ldots,k}=\sum_{0\ls s<r}\prod_{j=1}^{p-1}\bi{j(pm+s)}{pm+s}
\eq S\ (\mo\ p),$$
where $$S:=\sum_{0\ls s<\min\{r,2\}}\prod_{j=1}^{p-1}\(\bi{js}s\bi{jm}m\).$$
Clearly $S=0$ when $r=0$. If $r\gs 2$, then
$$S=(1+(p-1)!)\prod_{j=1}^{p-1}\bi{jm}m\eq0\ (\mo\ p).$$
In the case $r=1$ (i.e., $n\eq1\ (\mo\ p)$), if all the digits in the representation of
$n=pm+1$ in base $p$ belong to $\{0,1\}$, then
$$S=\prod_{j=1}^{p-1}\bi{jm}m\eq (-1)^{\psi_p(n)-1}\ (\mo\ p)$$
by Lemma 4.1, otherwise $S\eq0\ (\mo\ p)$ in view of Lemma 4.1. This ends the proof of part (ii).

(iii) Part (iii) of Theorem 1.2 follows immediately from Lemma 2.4.

 By the above we have completed the proof of Theorem 1.2. \qed

\medskip

\Ack. The author is grateful to the referee for many helpful comments.

\enddocument